\documentclass[letterpaper, 10 pt, conference]{ieeeconf}  
\IEEEoverridecommandlockouts                              
\usepackage{mathrsfs} 
\usepackage{kpfonts,comment} 
\usepackage{times}   
\usepackage{cancel}
\usepackage{kz_style}
\usepackage{setspace}
\usepackage{booktabs}
\usepackage{tabularx} 
\usepackage{siunitx}
\usepackage{sidecap}
\usepackage{{makecell}}
\setcellgapes{3pt}\makegapedcells
\usepackage{cite}
\renewcommand\theequation{\arabic{section}.\arabic{equation}}

\renewcommand\thesection{\arabic{section}}
\title{\LARGE \bf
Learning the Kalman Filter with Fine-Grained Sample Complexity}

\author{Xiangyuan Zhang \qquad Bin Hu \qquad Tamer Ba\c{s}ar
\thanks{The authors are with the Department of Electrical and Computer Engineering and Coordinated Science Laboratory, University of Illinois at Urbana-Champaign, Urbana, IL 61801, USA
        {\tt\small \{xz7, binhu7, basar1\}@illinois.edu}. X. Zhang and T. Ba\c{s}ar were supported in part by the US Army Research Laboratory (ARL) Cooperative Agreement W911NF-17-2-0196, and in part by the Office of Naval Research (ONR) MURI Grant N00014-16-1-2710. B. Hu was supported by the National Science Foundation (NSF) award CAREER-2048168. }%
}

\begin{document}

\maketitle
\thispagestyle{empty}
\pagestyle{empty}

\begin{abstract}
We develop the first end-to-end sample complexity of model-free policy gradient (PG) methods in discrete-time infinite-horizon Kalman filtering. Specifically, we introduce the receding-horizon policy gradient (RHPG-KF) framework and demonstrate $\tilde{\mathcal{O}}(\epsilon^{-2})$ sample complexity for RHPG-KF in learning a stabilizing filter that is $\epsilon$-close to the optimal Kalman filter. Notably, the proposed RHPG-KF framework does not require the system to be open-loop stable nor assume any prior knowledge of a stabilizing filter. Our results shed light on applying model-free PG methods to control a linear dynamical system where the state measurements could be corrupted by statistical noises and other (possibly adversarial) disturbances.
\end{abstract}

\section{Introduction}

In recent years, policy-based reinforcement learning (RL) methods \cite{sutton2000policy, kakade2002natural,schulman2015trust,schulman2017proximal} have gained increasing attention in continuous control applications \cite{schulman2015high,lillicrap2015continuous,recht2019tour}. While traditional model-based techniques synthesize controller designs in a case-by-case manner \cite{anderson1979optimal, anderson1990optimal}, model-free policy gradient (PG) methods promise a universal framework that learns controller designs in an end-to-end fashion. The universality of model-free PG methods makes them desired candidates in complex control applications that involve nonlinear system dynamics and imperfect state measurements. Despite countless empirical successes, the theoretical properties of model-free PG methods have not yet been fully investigated in continuous control. Initiated by \cite{fazel2018global}, a recent line of research has well-analyzed the sample complexity of zeroth-order PG methods in a number of linear state-feedback control benchmarks, including linear-quadratic regulator (LQR) \cite{mohammadi2019convergence, malik2020derivative, hambly2020policy, perdomo2021stabilizing, zhang2023revisit}, distributed/decentralized LQR \cite{li2019distributed, furieri2019learning}, and linear robust control \cite{gravell2019learning, zhang2021derivative, zhang2022rhpg}. However, the theoretical properties of PG methods remain elusive in the output-feedback control settings, where the state measurement process could be corrupted by statistical noises and/or other (possibly adversarial) disturbances. 

In this work, we take an initial step and study the sample complexity of zeroth-order PG methods in the discrete-time infinite-horizon Kalman filtering (KF) problem \cite{kalman1960new, anderson1979optimal}. The KF problem aims to design an optimal filter that generates estimates of the unknown system states over time, by utilizing a sequence of observed measurements corrupted by statistical noises. The KF problem has been recognized as one of the cornerstones of modern control theory \cite{basar2001control}. Furthermore, in the linear-quadratic Gaussian (LQG) problem, the separation principle \cite{astrom1971introduction} states that the optimal control law combines KF and LQR. Thus, KF is a fundamental benchmark for studying the sample complexity of model-free PG methods beyond state-feedback settings. 

Despite being the dual problem to noise-less LQR \cite{kalman1960general, astrom1971introduction}, the KF problem possesses a substantially more complicated optimization landscape from the model-free PG perspective, since KF itself is a dynamical system rather than a static matrix. Specifically, the optimization problem over dynamic filters might admit multiple suboptimal stationary points, and the optimal KF possesses a set of equivalent realizations up to similarity transformations \cite{zheng2021analysis, umenberger2022globally}. None of the above challenges appear when using model-free PG to learn a static LQR policy \cite{fazel2018global, mohammadi2019convergence, malik2020derivative, hambly2020policy}. As a result of the challenging landscape the filtering problem presents, only a few papers have focused on dynamic output-feedback settings. In particular, \cite{zheng2021analysis} has analyzed the optimization landscape of LQG and \cite{hu2022connectivity} has studied the optimization landscape of LQG with an additional $\cH_{\infty}$-robustness constraint. The most relevant work \cite{umenberger2022globally} has shown that an ``informativity-regularized'' PG method provably converges to an optimal dynamic filter in the continuous-time KF problem, assuming that the model is known. However, it is not clear if one can directly apply the techniques in \cite{umenberger2022globally} to the model-free setting and obtain any sample complexity guarantees. Moreover, \cite{umenberger2022globally} has assumed an open-loop stable system and prior knowledge to an initial filter satisfying the ``informativity'' condition, which is a more stringent condition than the stability of the closed loop. Thus, obtaining sample complexity of model-free PG methods in the KF problem has remained as a major challenge.

In this work, we introduce the receding-horizon PG framework (RHPG-KF) and establish $\tilde{\mathcal{O}}(\epsilon^{-2})$ sample complexity for the convergence of RHPG-KF to an $\epsilon$-optimal dynamic KF. By unifying model-free PG methods with dynamic programming (DP), the proposed RHPG-KF framework first decomposes the KF problem into a sequence of unconstrained strongly-convex subproblems. Then, the RHPG-KF framework solves these subproblems recursively using model-free PG methods. Finally, we demonstrate that the accumulated optimization errors from solving each of the subproblems can be carefully controlled, and as a result, the learned filter is $\epsilon$-close to the optimal KF. 

The RHPG-KF framework in our current work, combined with the RHPG framework for linear control in \cite{zhang2022rhpg, zhang2023revisit}, demonstrates the significant utilization of DP in streamlining the analyses of model-free PG methods in linear control and estimation. Due to the separation principle \cite{astrom1971introduction}, our results shed light on applying model-free PG methods in solving the LQG problem through a sequential design of control and estimation.

\subsection{Notations}\label{sec:notations}
For a square matrix $X$, we denote its trace, spectral norm, condition number, and spectral radius as $\Tr(X)$, $\|X\|$, $\kappa_X$, and $\rho(X)$ respectively. We define the $W$-induced norm of $X$ as
\begin{align*}
	\|X\|^2_W := \max_{z \neq 0} \frac{z^{\top}X^{\top}WXz}{z^{\top}Wz}.
\end{align*}
If $X$ is further symmetric,  we use $X > 0$, $X\geq 0$, $X\leq 0$, and $X <0$ to denote that $X$ is positive definite (pd), positive semi-definite (psd), negative semi-definite (nsd), and negative definite (nd), respectively. We use $x\sim \cN(\mu,\Sigma)$ to denote a Gaussian random vector with mean $\mu$ and covariance $\Sigma$. Lastly, we use $\bI$ and $\bm{0}$ to denote the identity and zero matrices, respectively, with appropriate dimensions.

\section{Preliminaries}
\subsection{Infinite-Horizon Kalman Filtering}
Consider the discrete-time linear time-invariant system
\begin{align}\label{eqn:inf_LQE_dynamics}
	x_{t+1} = Ax_t+ w_t, \quad y_t = Cx_t + v_t,
\end{align}
where $x_t \in \RR^n$ is the state, $y_t \in \RR^m$ is the output measurement, and $w_t \sim \cN(\bm{0}, W)$, $v_t\sim \cN(\bm{0}, V)$ are sequences of i.i.d. zero-mean Gaussian noises for some $W, V > 0$, also independent of each other. The initial state is also assumed to be a Gaussian random vector such that $x_0 \sim \cN(\bar{x}_0, X_0)$, independent of $\{w_t, v_t\}$, and $X_0 > 0$. Additionally, we assume that $(C, A)$ is observable and note that the condition $W>0$ readily leads to the controllability of $(A, W^{1/2})$, which is a standard condition in KF. 

The KF problem aims to generate a sequence of estimated states, denoted by $\hat{x}_t$ for each $t$, that minimizes the infinite-horizon mean-square error (MSE):
\begin{align}\label{eqn:inf_mse}
	\cJ_{\infty} := \lim_{N\to \infty}\frac{1}{N}\EE\bigg\{\sum_{t=0}^{N}(x_t-\hat{x}_t)^{\top}(x_t-\hat{x}_t)\bigg\}.
\end{align}
Moreover, each $\hat{x}_t$ can only depend on the history and output measurements up to but not including $t$, i.e., $\{y_0, \cdots, y_{t-1}\}$. The celebrated result of Kalman \cite{kalman1960new} showed that the $\cJ_{\infty}$-minimizing filter (could also be called $1$-step predictor), which exists under the controllability and the observability conditions, has the form of 
\begin{align}\label{eqn:dynamic_filter}
	\hat{x}_{t+1}^* &= (A-L^*C)\hat{x}_t^* + L^*y_t, \quad  \hat{x}_{0}^* = \bar{x}_0, \\ 
	L^*&=A\Sigma^* C^{\top}(V + C\Sigma^* C^{\top})^{-1}, \label{eqn:kalman_gain}
\end{align} 
where $L^*$ is the Kalman gain and $\Sigma^*$ represents the unique pd solution to the filter algebraic Riccati equation (FARE):
\begin{align}
	\Sigma &= A\Sigma A^{\top} - A\Sigma C^{\top}(V + C\Sigma C^{\top})^{-1}C\Sigma A^{\top} + W. \label{eqn:filter_riccati}
\end{align}
Hence, without any loss of optimality, we can restrict the search to the class of filters of the form $\hat{x}_{t+1} = A_L\hat{x}_t + B_Ly_t$ and then parametrize the KF problem as a minimization problem over $A_L$ and $B_L$ subject to a stability constraint\footnote{Extending the results in this work to the setting with instantaneous feedback measurement (i.e., allowing $\hat{x}_t$ to depend also on $y_t$, and hence replacing $y_t$ in (\ref{eqn:filter_form}) with $y_{t+1}$) would be straightforward.}
\begin{align}\label{eqn:filter_form}
	&\min_{A_L, B_L} \quad \cJ_{\infty}(A_L, B_L) \\
 \text{s.t.} \quad &\hat{x}_{t+1} = A_L\hat{x}_t + B_Ly_t ~\ \text{and} ~\ \rho(A_L) < 1. \nonumber
\end{align}
Note that by \eqref{eqn:dynamic_filter}, there indeed exists a solution to \eqref{eqn:filter_form} where $(A_L^*, B_L^*)=(A-L^*C, L^*)$. Note also that when $(A, C)$ is known, \eqref{eqn:filter_form} involves an over-parametrization since solving \eqref{eqn:filter_form} is equivalent to optimizing a single variable $B_L$. However, in the model-free setting where $(A, C)$ is unknown, which is the target setting of our paper, it is reasonable to parametrize the KF problem as in \eqref{eqn:filter_form}. Until now, obtaining sample complexity of model-free PG methods in solving the KF problem \eqref{eqn:filter_form} has remained as a major challenge.

\begin{figure*}[t]
\centering 
\includegraphics[width = 0.8\textwidth]{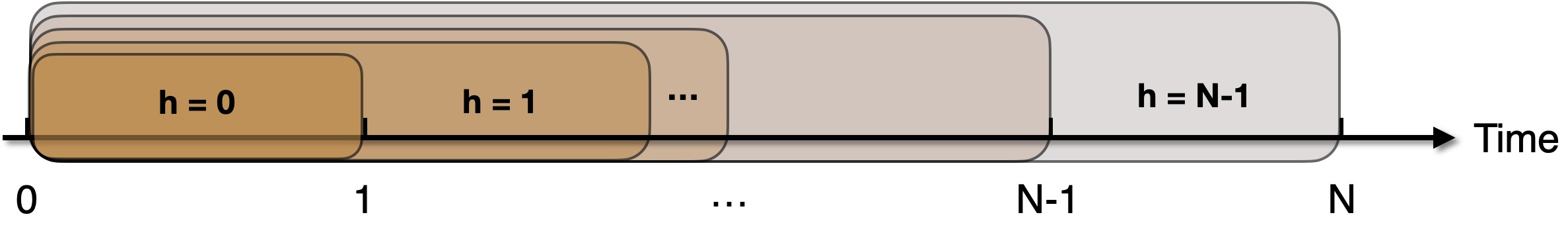}	
\caption{Illustrating the procedure of Algorithm \ref{alg:DPfilter}, which solves, at each iteration indexed by $h$, an $(h+1)$-horizon KF problem from $t=0$ to $t=h+1$, where we fix the filter parameters for all $t\in\{0, \cdots, h-1\}$ as those generated from earlier iterations and only optimize for filter parameters $(A_{L_h}, B_{L_h})$. }\label{fig:DP_illustration}
\end{figure*}

\subsection{Finite-Horizon Kalman Filtering}
We present the finite-$N$-horizon version of the KF problem, which is also described by the system dynamics \eqref{eqn:inf_LQE_dynamics}. Adopting the same parametrization as in \eqref{eqn:filter_form}, but this time allowing time-dependence, and again without any loss of optimality, we represent the finite-horizon KF problem as a minimization problem over a sequence of time-varying filter parameters $\{A_{L_t}, B_{L_t}\}$, for all $t\in\{0, \cdots, N-1\}$, 
\begin{align}\label{eqn:finite-mse}
	\min_{\{A_{L_t}, B_{L_t}\}} \ &\cJ\big(\{A_{L_t}, B_{L_t}\}\big) := \EE\bigg\{\sum_{t=0}^{N}(x_t-\hat{x}_t)^{\top}(x_t-\hat{x}_t)\bigg\} \\
	&\text{s.t.} \quad \hat{x}_{t+1} = A_{L_t}\hat{x}_t + B_{L_t}y_t, \quad \hat{x}_0 = \bar{x}_0. \nonumber
\end{align}
The minimum in \eqref{eqn:finite-mse} can be achieved by $(A^*_{L_t}, B^*_{L_t}) = (A-L^*_tC, L^*_t)$, where $L^*_t$ is the time-varying Kalman gain
\begin{align}
	&L_t^*=A\Sigma_t^* C^{\top}(V + C\Sigma_t^* C^{\top})^{-1}, \quad  \Sigma_0^* = X_0.\label{eqn:kalman_gain_finite}\\
	\hspace{-0.8em} \Sigma^*_{t+1} &= A\Sigma^*_{t}A^{\top} - A\Sigma^*_t C^{\top}(V \hspace{-0.1em}+\hspace{-0.1em} C\Sigma^*_t C^{\top})^{-1}C\Sigma^*_t A^{\top} \hspace{-0.1em}+\hspace{-0.1em} W. \label{eqn:filter_riccati_finite}
\end{align}
The solutions $\Sigma^*_t$, for all $t \in \{0, \cdots, N-1\}$, generated by the filter Riccati difference equation (FRDE) \eqref{eqn:filter_riccati_finite} always exist and are unique and pd, due to $V > 0$, $W > 0$, and the iteration starting with $\Sigma^*_0 = X_0 > 0$.

\section{Receding-Horizon Policy Gradient}
\subsection{Kalman Filtering with Dynamic Programming}
It is well known that the solution of the FRDE \eqref{eqn:filter_riccati_finite} converges monotonically to the stabilizing solution of the FARE \eqref{eqn:filter_riccati}, and the rate of this convergence is exponential \cite{chan1984convergence, hassibi1999indefinite}. Then, it readily follows that the optimal time-varying filter parameters $(A_{L_t}^*, B_{L_t}^*)$ to the finite-horizon KF problem \eqref{eqn:finite-mse} also converge monotonically to the time-invariant $(A_{L}^*, B_{L}^*)$ as $N \to \infty$. Before formally presenting the convergence result, we assume that $\bar{x}_0 \neq 0$ and $X_0 > \Sigma$, which could be fulfilled in practice by injecting an additional Gaussian noise vector into $x_0$. This additional noise vector will not affect the solution to the infinite-horizon KF problem \eqref{eqn:kalman_gain}-\eqref{eqn:filter_riccati} since it is independent of $\bar{x}_0$ and $X_0$. Then, we characterize the exponential rate of the convergence in policy (i.e., $(A_{L_t}^*, B_{L_t}^*) \to (A_L^*, B_L^*)$) in the following theorem.

\begin{theorem}\label{lemma:finite_approximation}
The finite-horizon Kalman gain as in \eqref{eqn:kalman_gain_finite} converges  to the infinite-horizon Kalman gain defined in \eqref{eqn:kalman_gain} exponentially fast as $N\to \infty$. Specifically, using $\|\cdot\|_*$ to denote the $\Sigma^*$-induced norm and letting
	\begin{align}\label{eqn:N0}
		N_0 = \frac{1}{2}\cdot \frac{\log\big(\frac{\|X_0-\Sigma^*\|_*\cdot\kappa_{\Sigma^*}\cdot \|A_L^*\|\cdot\|C\|} {\epsilon\cdot\lambda_{\min}(V)}\big)}{\log\big(\frac{1}{\|A_L^*\|_*}\big)} + 1.
	\end{align}
	where $\|A_L^*\|_* <1$, we have that, for all $N\geq N_0$, it holds that $L^*_{N-1}$ is stabilizing and $\|L^*_{N-1} - L^*\| \leq \epsilon$ for any $\epsilon > 0$.
\end{theorem}\par  
The proof of Theorem \ref{lemma:finite_approximation} is provided in \S\ref{proof:finite}. Theorem \ref{lemma:finite_approximation} demonstrates that if  $N\sim \cO(\log(\epsilon^{-1}))$, then solving the finite-horizon KF problem will return filter parameters $A_{L_{N-1}}^*$ and $B_{L_{N-1}}^*$ that are $\epsilon$-close to $A_L^*$ and $B_L^*$, respectively. Furthermore, it holds that $\rho(A_{L_{N-1}}^*)<1$ for any $\epsilon > 0$.

\begin{algorithm}[t]
  \label{alg:DPfilter}
 \caption{RHPG-KF}
  \SetAlgoLined
  \KwIn{horizon $N$, max iterations $\{T_0, \cdots, T_{N-1}\}$}
  \For{$h = 0, \cdots, N-1$}{ 
Initialize $[A_{{L_h}, 0} \ B_{{L_h}, 0}] \leftarrow [\bm{0} \ \bm{0}]$\;

\For{$i = 0, \cdots, T_h-1$}{
    Execute \eqref{eqn:PG} to compute $[A_{L_{h}, i+1} \ B_{L_{h}, i+1}]$\;
  }
  }
  Return $A_{L_{N-1}, T_{N-1}}, B_{L_{N-1}, T_{N-1}}$\;
\end{algorithm}

\subsection{Algorithm Design}
Instead of solving \eqref{eqn:finite-mse} directly, we propose Algorithm \ref{alg:DPfilter} that sequentially solves $N$ optimization problems, as illustrated in Figure \ref{fig:DP_illustration}. In particular, for every $h \in \{0, \cdots, N-1\}$, Algorithm \ref{alg:DPfilter} solves a KF problem, denoted as $\cF_h$, that starts at time $t=0$ with a randomly sampled $x_{0} \sim \cN(\bar{x}_0, X_0)$ and ends at $t=h+1$. Moreover, in solving every $\cF_h$ we only optimize for the one-step filter parameters $(A_{L_h}, B_{L_h})$, while fixing the prior filter parameters $\{(A_{L_t}, B_{L_t})\}$ for all $t\in\{0, \cdots, h-1\}$ to be the convergent ones from earlier iterations of Algorithm \ref{alg:DPfilter}. This procedure renders each iteration of Algorithm \ref{alg:DPfilter} into a \emph{quadratic} program in $(A_{L_h}, B_{L_h})$. 

The intuition of Algorithm \ref{alg:DPfilter} comes from Bellman's principle of optimality. More explicitly, for every $h\in\{0, \cdots, N-1\}$, the truncation of the entire sequence of optimal filter parameters $\{(A_{L_t}^*, B_{L_t}^*)\}_{t\in\{0, \cdots, N-1\}}$ to the interval $[0, h]$ constitutes the optimal filter of the KF problem $\cF_h$. Based on Bellman's principle of optimality, we can decompose \eqref{eqn:finite-mse} into sequentially solving a series of (one-step) KF problems $\{\cF_h\}_{h\in\{0, \cdots, N-1\}}$, assuming (until \S\ref{sec:bias}) that each one of them can be solved exactly.

Concretely, for every iteration $h\in\{0, \cdots, N-1\}$, Algorithm \ref{alg:DPfilter} solves a one-step minimization problem
\begin{align}
	\min_{A_{L_h}, B_{L_h}} \ &\cJ_h(A_{L_h}, B_{L_h}) := \EE\bigg[\sum_{t=0}^{h+1}(x_t-\hat{x}_t)^{\top}(x_t-\hat{x}_t)\bigg] \label{eqn:induction}\\
	&\text{s.t.} \quad \hat{x}_{h+1}= A_{L_h}\hat{x}_h + B_{L_h}y_h.\label{eqn:induction_cons}
\end{align}
where $\hat{x}_h$ in \eqref{eqn:induction_cons} is generated by applying the convergent filter parameters in previous iterations of Algorithm \ref{alg:DPfilter}, i.e., $\{(A_{L_t}^*, B_{L_t}^*)\}$ for all $t \in \{0, \cdots, h-1\}$ and thus is independent of $(A_{L_h}, B_{L_h})$. Moreover, $y_h$ in \eqref{eqn:induction_cons} is computed by \eqref{eqn:inf_LQE_dynamics} and therefore is also independent of $(A_{L_h}, B_{L_h})$. As a result, \eqref{eqn:induction} is a quadratic unconstrained minimization problem in $(A_{L_h}, B_{L_h})$. This implies that one can simply apply any PG methods with the initial point being $(\bm{0}, \bm{0})$ to solve \eqref{eqn:induction}. Lastly, to solve \eqref{eqn:induction} using model-free PG methods for all $h$, it is standard to impose the following assumption.
\begin{assumption}\label{assumption:simulator}
The user has access to a simulator (e.g., Algorithm \ref{alg:estimator}) such that for any input filter policies $(A_{L_h}, B_{L_h})$, the simulator can return an empirical objective value \eqref{eqn:induction}.	
\end{assumption}

Assumption \ref{assumption:simulator} requires the simulator to generate exact state trajectories of the simulated model, but it only reveals a noisy scalar objective value to the learning algorithm. This assumption is reasonable for model-free learning since the algorithm does not use any system information directly. However, building a simulator naturally requires some knowledge of the system model, which could be either exact, or approximate, or simplified. The focus of our paper is to investigate the theoretical properties of model-free PG methods in learning the KF with Assumption \ref{assumption:simulator} standing. Whether and how Assumption \ref{assumption:simulator} could be relaxed is left as an important future research topic.

\subsection{Bias of Model-Free Receding-Horizon Filtering}\label{sec:bias}
The procedure of Algorithm \ref{alg:DPfilter} relies on Bellman's principle of optimality, which requires solving each of the iterations exactly. However, iterative algorithms such as PG methods can only return an $\epsilon$-accurate solution in a finite time. Hence, it is paramount to analyze how computational errors accumulate in the forward DP process. Accumulation of errors also brings up the question of the accuracy of the output filter from the exact solution $(A_L^*, B_L^*)$ to the infinite horizon KF problem. In the theorem below, we show that if we choose $N\sim \cO(\log(\epsilon^{-1}))$ by Theorem \ref{lemma:finite_approximation} and apply a PG method that solves each of the one-step KF problems to an $\epsilon$-accurate solution, then the output filter of Algorithm \ref{alg:DPfilter} would be $\epsilon$-close to $(A_L^*, B_L^*)$ of the infinite-horizon problem. \par 

\begin{theorem}\label{theorem:KF_DP}
	Choose $N$ according to Theorem \ref{lemma:finite_approximation} and assume that one can compute, for all $h\in\{0, \cdots, N-1\}$ and some $\epsilon > 0$, a pair of matrices $(\tilde{A}_{L_h}, \tilde{B}_{L_h})$ that satisfies 
	\begin{align*}
		\big\|\tilde{A}_{L_h}  - \tilde{A}_{L_h}^*\big\|, \big\|\tilde{B}_{L_h}  - \tilde{B}_{L_h}^*\big\| \hspace{-0.1em}\sim\hspace{-0.1em}\cO(\epsilon\cdot \texttt{poly}(\text{system parameters})),
	\end{align*}
	where $(\tilde{A}_{L_h}^*, \tilde{B}_{L_h}^*)$ constitutes the optimal filter of $\cF_{h}$ after absorbing errors in all previous iterations of Algorithm \ref{alg:DPfilter}. Then, Algorithm \ref{alg:DPfilter} outputs a pair $(\tilde{A}_{L_{N-1}}, \tilde{B}_{L_{N-1}})$ that satisfies $\big\|[\tilde{A}_{L_{N-1}} \ \tilde{B}_{L_{N-1}}] - [A^*_L \ B^*_L]\big\| \leq \epsilon$. Furthermore, if $\epsilon$ is sufficiently small such that $\epsilon < 1-\|A_L^*\|_*$, then $\tilde{A}_{L_{N-1}}$ is guaranteed to be stabilizing.
\end{theorem}

\begin{figure*}
	\centering
	\includegraphics[width=0.98\textwidth]{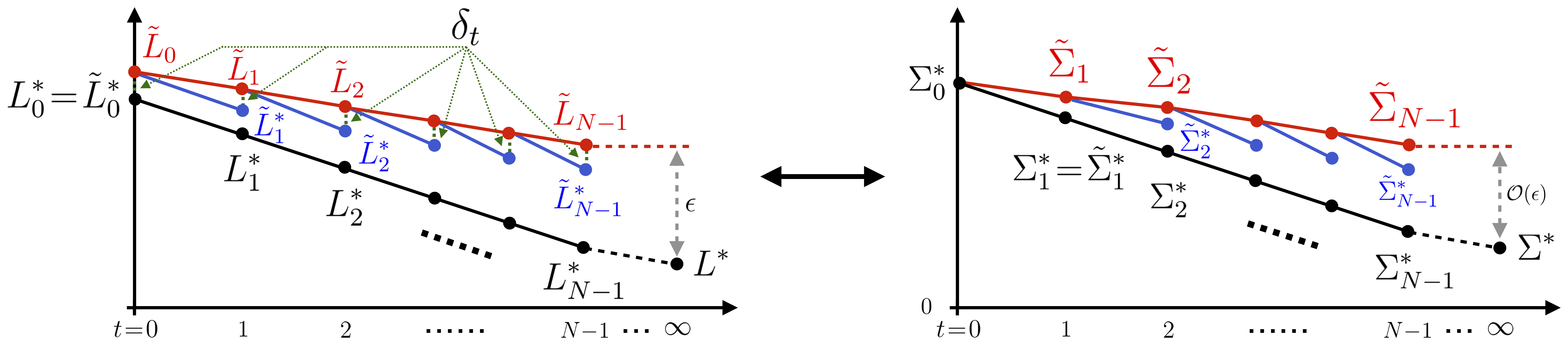}
	\caption{Illustrating the idea of Theorem \ref{theorem:KF_DP} in the scalar case, where we aim to obtain a filter gain $\tilde{L}_{N-1}$ that is $\epsilon$-close to $L^*$. Note that the convergence of $\tilde{L}_{N-1}$ to $L^*$ readily implies the convergence of $(\tilde{A}_{L_{N-1}}, \tilde{B}_{L_{N-1}})$ to $(A^*_L, B_L^*)$. First, Theorem \ref{lemma:finite_approximation} proves that $L^*_{N-1}$ is $\epsilon$-close to $L^*$ by selecting $N$ accordingly. Subsequently, Theorem \ref{theorem:KF_DP} analyzes the propagation of the optimization errors (forward in time) from solving each of the one-step subproblems, denoted as $\delta_t:=\tilde{L}_t - \tilde{L}^*_t$ for all $t$, where $\tilde{L}^*_t$ represents the current optimal filter gain after absorbing errors from all previous iterations. We demonstrate that if $\|\delta_t\|\sim\cO(\epsilon)$ for all $t$, then the convergent filter gain $\tilde{L}_{N-1}$ returned by Algorithm \ref{theorem:KF_DP} is guaranteed to be $\epsilon$-close to the optimal KF gain $L^*$.}\label{fig:proof_sketch}
	\vspace{-0.4em}
\end{figure*}

We illustrate the idea of Theorem \ref{theorem:KF_DP} in Figure \ref{fig:proof_sketch} and defer its proof to \S\ref{proof:KF_DP}. Theorem \ref{theorem:KF_DP} guarantees that if every iteration of the DP is solved to an $\cO(\epsilon)$-accuracy, then the convergent filter after completing the DP procedure is at most $\epsilon$-away from the exact solution to the infinite-horizon KF. Then, it remains to show that (zeroth-order) PG methods do converge in every iterations of Algorithm \ref{alg:DPfilter}.

We note that Algorithm \ref{alg:DPfilter} takes two layers of approximation in solving LQ estimation problems (cf., Figure \ref{fig:proof_sketch}). First, we approximate the solution to the infinite-horizon KF problem \eqref{eqn:filter_form} by the solution to a finite-horizon simplification \eqref{eqn:finite-mse}, where we choose $N\sim \cO(\log(\epsilon^{-1}))$ following the exponential attraction of the Riccati equation. Then, we solve the finite-horizon problem \eqref{eqn:finite-mse} by combining forward DP with an arbitrary model-free solver. As a result of these two layers of approximation, our proposed framework appears to be the first that solves the infinite-horizon KF problem using only oracle-level accesses to the system.

\begin{algorithm}[t]
\DontPrintSemicolon
\caption{~Two-Point Zeroth-Order Oracle} 
  \SetAlgoLined
	\label{alg:estimator}
	\KwIn{$\bL_h:= [A_{L_h} \ B_{L_h}]$, smoothing radius $r_h$}
	1. Sample $\bL_h^{+} = \bL_h + r_hU$ and $\bL_h^{-} = \bL_h - r_hU$, where $U$ is uniformly sampled from the surface of an $n(m+n)$-dimensional unit sphere, i.e., $\|U\|_F=1$. \;
	2. Sample a random $x_0\sim\cN(\bar{x}_0, X_0)$, apply the filter with parameters being $\{\bL_t\}_{t\in\{0, \cdots, h-1\}}$, and arrive at state $x_h$ with prediction $\hat{x}_{h}$. \;
	3. Sample $w_h \sim \cN(\bm{0}, W)$ and $v_h \sim \cN(\bm{0}, V)$ and generate $x_{h+1}$ and $y_h$. \;
	4. Simulate two one-step filters with parameters being $\bL_h^{+}$ and $\bL_h^{-}$ and compute the corresponding estimates $\hat{x}_{h+1}^{+}$ and $\hat{x}_{h+1}^{-}$. Calculate the empirical MSEs $J_h(\bL_h^{+})$ and $J_h(\bL_h^{-})$.\;
	5. Return: {\small$g_h = \frac{n(m+n)}{2r_h}\cdot\big[J_h(\bL_h^{+})-J_h(\bL_h^{-})\big]\cdot U$.}\;
\end{algorithm}

\section{Convergence and Sample Complexity}
We are now ready to present the sample complexity of Algorithm \ref{alg:DPfilter} by analyzing the number of samples required in solving each iteration. First, we summarize the optimization landscapes of each quadratic program in Algorithm \ref{alg:DPfilter}.   
\begin{fact}\label{lemma:convex}
There exist constants $m, \psi > 0$ such that for all $h\in\{0, \cdots, N-1\}$, the objective function $\cJ_h(A_{L_h}, B_{L_h})$ in \eqref{eqn:induction} is $m$-strongly convex and $\psi$-smooth in $(A_{L_h}, B_{L_h})$.
\end{fact}
Fact \ref{lemma:convex} is a direct result of the quadratic objective function $\cJ_h((A_{L_h}, B_{L_h}))$. Then, we define the vanilla PG update as\footnote{Our framework is also compatible with other PG methods such as natural PG \cite{kakade2002natural} and least-squares policy iteration \cite{lagoudakis2003least}.}
\small
\begin{align}
	 \big[A_{L_h}' \ B_{L_h}'\big] =  \big[A_{L_h} \ B_{L_h}\big] - \eta_h \cdot  \Big[\frac{\partial \cJ_h(A_{L_h}, B_{L_h})}{\partial A_{L_h}} \ \frac{\partial \cJ_h(A_{L_h}, B_{L_h})}{\partial B_{L_h}}\Big], \label{eqn:PG}
\end{align}
\normalsize
where $\eta_h > 0$ is a constant stepsize to be determined. We first provide the convergence analysis of \eqref{eqn:PG} in the exact case (i.e., when exact gradients are accessible) based on standard results in minimization over strongly-convex and smooth functions.

\begin{fact}\label{lemma:convergence}
For all $h \in \{0, \cdots, N-1\}$ and a fixed $\epsilon > 0$, choose a constant stepsize $\eta_h \in (0,  1/\psi]$. Then, the deterministic PG update \eqref{eqn:PG} converges linearly to $[\tilde{A}_{L_h}^* \ \tilde{B}^*_{L_h}]$ such that $\big\|[A_{{L_h}, T_h} \ B_{{L_h}, T_h}] - [\tilde{A}_{L_h}^* \ \tilde{B}^*_{L_h}]\big\| \leq \epsilon$ after a total number of $T_h\sim\cO(\log(\epsilon^{-1}))$ iterations.
\end{fact}

When the exact PGs are not available, \eqref{eqn:PG} can be implemented using estimated PGs sampled from system trajectories using (two-point) zeroth-order optimization techniques (cf., Algorithm \ref{alg:estimator}). Notably, the convergence rate of zeroth-order PG update is at most $\cO(\sqrt{n(m+n)})$ slower than that of the \emph{stochastic} first-order PG update \cite{duchi2015optimal}, where $n(m+n)$ is the dimension of $[A_{L_h} \ B_{L_h}]$ for all $h$. Furthermore, in strongly convex and smooth minimization problems, the stochastic PG update converges to an $\epsilon$-optimal policy with a high probability at the rate of $\cO(\epsilon^{-2})$ (e.g., Proposition 1 of \cite{rakhlin2011making}). We combine these two results to obtain the sample complexity bound in the following proposition.

\begin{proposition}\label{prop:sample}
	For all $h\in\{0, \cdots, N-1\}$, choose the smoothing radius $r_{h, i} \sim \cO(\sqrt{\epsilon}\cdot i^{-1})$ and the stepsize $\eta_{h, i}\sim \cO(i^{-1})$, where $i$ is the iteration index. Then, the zeroth-order PG update \eqref{eqn:PG} converges after $T_h \sim \tilde{\cO}(\epsilon^{-2}\cdot\log(\delta^{-1}))$ iterations in the sense that $\big\|[A_{{L_h}, T_h} \ B_{{L_h}, T_h}] - [\tilde{A}_{L_h}^* \ \tilde{B}^*_{L_h}]\big\| \leq \epsilon$ with a probability of at least $1-\delta$.
\end{proposition}

Combining Theorem \ref{theorem:KF_DP} with Proposition \ref{prop:sample}, we conclude that if we spend $\tilde{O}(\epsilon^{-2})$ samples in solving every one-step KF problem to $\cO(\epsilon)$-accuracy with a probability of $1-\delta$, for all $h \in \{0, \cdots, N-1\}$, then Algorithm \ref{alg:DPfilter} is guaranteed to output a pair $(\tilde{A}_{L_{N-1}}, \tilde{B}_{L_{N-1}})$ that satisfies 
\begin{align*}
	\Big\|\big[\tilde{A}_{L_{N-1}} \ \tilde{B}_{L_{N-1}}\big] - \big[A_{L}^* \ B^*_{L}\big]\Big\| \leq \epsilon
\end{align*}
with a probability of at least $1-N\delta$.  The total sample complexity of Algorithm \ref{alg:DPfilter} is thus $\tilde{O}(\epsilon^{-2})\cdot O(\log(\epsilon^{-1}))\sim \tilde{O}(\epsilon^{-2})$. Lastly, we compare our sample complexity with those obtained in LQR \cite{malik2020derivative} in the following remark.
\begin{remark}
	 Our $\tilde{O}(\epsilon^{-2})$ sample complexity bound (for the convergence in policy) matches the best-known sample complexity of two-point zeroth-order PG methods in LQR \cite{malik2020derivative}. This is due to that the $\tilde{\cO}(\epsilon^{-1})$ sample complexity for the convergence in the objective value (e.g., $f(K) - f(K^*) \leq \epsilon$) in \cite{malik2020derivative} is equivalent to an $\tilde{\cO}(\epsilon^{-2})$ sample complexity for the convergence in policy (i.e., $\|K-K^*\| \leq \epsilon$). (cf., Sec. A.10 in page 33 of \cite{zhang2021derivative}). 
\end{remark}

\section{Numerical Experiments}
To complement our theoretical results, we perform simulations on a scalar infinite-horizon KF problem with:
\begin{align}\label{eqn:sim_kf}
	A = 2, \quad C =  W =  V = 1, \quad  \bar{x}_0 = 1, \quad X_0 = 5.
\end{align}
The discrete-time system \eqref{eqn:sim_kf} is open-loop unstable and the optimal filter is $(A_L^*, B_L^*) = (0.3820, 1.6180)$.  We choose $N = \texttt{ceil}(\log(\epsilon^{-1}))$ and manually select a constant stepsize in each iteration of Algorithm \ref{alg:DPfilter}. Furthermore, we choose $r_h = \sqrt{\epsilon}$ and run Algorithm \ref{alg:DPfilter} to solve the KF problem under six different required convergence accuracies $\epsilon \in \{10^{-3}, 3.16\times 10^{-3}, 10^{-2}, 3.16\times 10^{-2}, 10^{-1}, 3.16\times 10^{-1}\}$. Specifically, in every iteration $h$ of Algorithm \ref{alg:DPfilter}, we initialize the PG update with the filter parameters being $(A_{L_h}, B_{L_h}) = (0, 0)$ and apply the zeroth-order update until convergence in the sense that $\big\|[\tilde{A}_{{L_h}} \ \tilde{B}_{{L_h}}] - [\tilde{A}_{L_h}^* \ \tilde{B}^*_{L_h}]\big\| \leq \epsilon/N$.

\begin{figure}
	\centering\hspace{-1.5em}
	\includegraphics[width=0.46\textwidth]{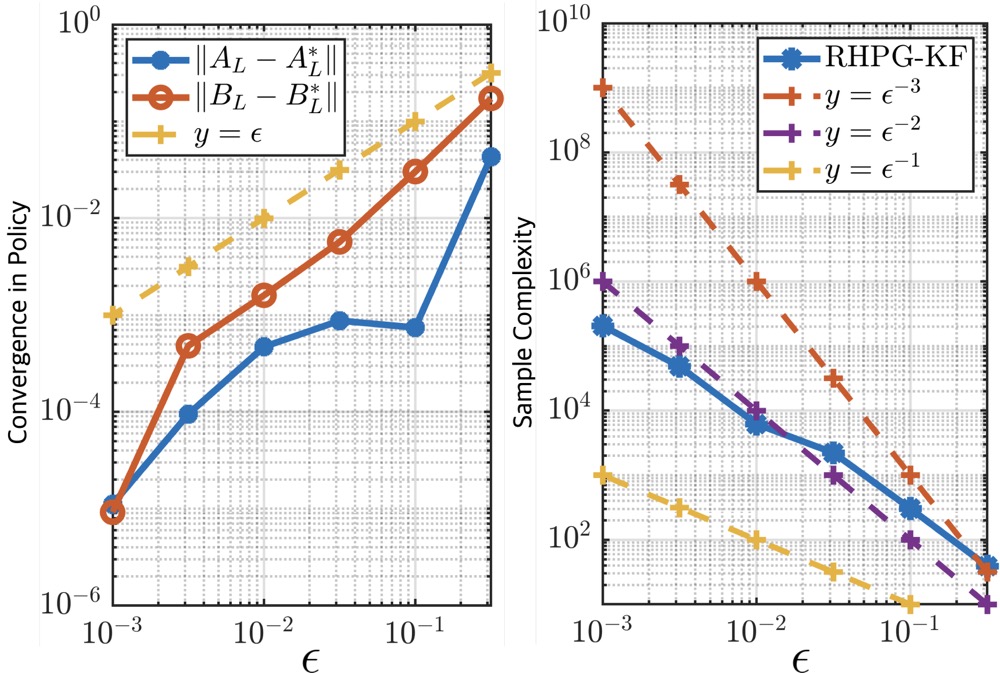}
	\caption{Left: error differences between the convergent filter parameters returned from Algorithm \ref{alg:DPfilter} and $A_L^*, B_L^*$. Right: the total number of samples (i.e., the number of calls to the zeroth-order oracle) used by Algorithm \ref{alg:DPfilter} under six different required accuracies.}\label{fig:KF_sim}
	\vspace{-1.5em}
\end{figure}

As shown in the left graph of Figure \ref{fig:KF_sim}, Algorithm \ref{alg:DPfilter} successfully returns a convergent filter that is $\epsilon$-close to the optimal KF, for all six different choices of $\epsilon$. Moreover, for all six cases, the convergent filter is stabilizing. Then, the right graph of Figure \ref{fig:KF_sim} verifies the $\tilde{\cO}(\epsilon^{-2})$ sample complexity bound presented in Proposition \ref{prop:sample}.

To further demonstrate the scalability of the RHPG-KF framework, we apply Algorithm \ref{alg:DPfilter} to a vector example
\small
\begin{align}
	&\hspace{2em}A = \begin{bmatrix}
		9.9 & -0.02\\
		0.01 & 10.1
	\end{bmatrix}, \quad C = \begin{bmatrix}
		0.99 & 0\\
		-0.01 & 1.01	\end{bmatrix}, \label{eqn:vector_example}\\
		&W = 10^{-3}\cdot\bI, \quad V = 10^{-2}\cdot\bI, \quad \bar{x}_0 = \begin{bmatrix}
			0.1 \\ 0.1
		\end{bmatrix}, \quad X_0 = 2\cdot\bI. \nonumber
\end{align}
\normalsize
It can be computed that $\rho(A) = 10.0990$. The optimal KF gain is $\small{L^* = \begin{bmatrix}
		9.8979 & -0.0197 \\ 0.1099 & 9.9021
	\end{bmatrix}}$. We set $\epsilon=0.8$, $r=0.01$, and select $N=2$ by \eqref{eqn:N0}. We also initialize the zeroth-order PG update \eqref{eqn:PG} with $(\bm{0}, \bm{0})$ for all $h\in \{0, 1\}$, and emphasize that this initialization could be arbitrary. We run the zeroth-order update until convergence in the sense that $\big\|[\tilde{A}_{{L_h}} \ \tilde{B}_{{L_h}}] - [\tilde{A}_{L_h}^* \ \tilde{B}^*_{L_h}]\big\| \leq \epsilon/2$ for all $h$. Algorithm \ref{alg:DPfilter} successfully returns a pair $(\tilde{A}_{L_1}, \tilde{B}_{L_1})$ that satisfies $\rho(\tilde{A}_{L_1}) = 0.0341<1$ and $\big\|[\tilde{A}_{{L_1}} \ \tilde{B}_{{L_1}}] - [A_{L}^* \ B^*_{L}]\big\| = 0.4067 < \epsilon$.

\section{Discussion and Future Directions}\label{sec:discuss}
The RHPG-KF framework unifies DP with model-free PG methods by sequentially decomposing the KF problems into a series of unconstrained quadratic optimization problems. Together with \cite{zhang2022rhpg, zhang2023revisit}, we have demonstrated the significant utilization of DP in streamlining/simplifying the analyses of model-free PG methods in linear control and estimation. Yet another advantage of the RHPG framework is that it requires no assumptions nonstandard to model-free learning (e.g., knowing a priori a stabilizing initial point). Due to the broad applicability of DP and the underlying Bellman's principle of optimality, the proposed framework has the potential to be generalized to a wide variety of domains without any major changes to the underlying process.

The RHPG-KF framework can serve as a building block for applying model-free PG methods to solve the LQG problem, due to the separation principle \cite{astrom1971introduction}. One key observation is that the RHPG-KF framework does not require the system to be open-loop stable. Thus, one can learn the (central) LQG solution by solving the estimation and control problems separately. A possible approach is to fix a zero control and utilize the RHPG-KF framework to learn a KF. Subsequently, one can apply the algorithm in \cite{zhang2023revisit} to compute a certainty-equivalent control policy. However, the extension to LQG will not be straightforward due to the following reasons. First, LQG does not have any guaranteed margins \cite{doyle1978guaranteed}, and thus there exist LQG instances where an arbitrarily small inaccuracy (e.g., optimization error) to the optimal solution will de-stabilize the system. Moreover, the separation principle only holds when both KF and LQR are solved to the exact optimum, and generalizing it to the approximate setting requires non-trivial work. 


\section{Conclusion}
In this work, we have investigated the sample complexity of model-free PG methods for learning the discrete-time infinite-horizon KF. Our results serve as an initial step toward understanding the theoretical properties of model-free PG methods in the control of dynamical systems under uncertainty and using imperfect state measurements, which include the LQG problem as a prominent example.


\typeout{}
\small
\bibliographystyle{alpha}
\bibliography{main}
\normalsize

\renewcommand{\theequation}{\thesubsection.\arabic{equation}}

\appendix
\subsection{Proof of Theorem \ref{lemma:finite_approximation}}\label{proof:finite}
We first present a technical lemma due to \cite{de1989monotonicity}.
\begin{lemma} \label{lemma:comparison}
	Consider two RDEs 
\begin{align*}
	\Sigma_{t+1}^1 = W + A\Sigma_{t}^1A^{\top} - A\Sigma_{t}^1C^{\top}(V + C\Sigma_{t}^1C^{\top})^{-1}C\Sigma_{t}^1A^{\top}, \\
	\Sigma_{t+1}^2 = W + A\Sigma_{t}^2A^{\top} - A\Sigma_{t}^2C^{\top}(V + C\Sigma_{t}^2C^{\top})^{-1}C\Sigma_{t}^2A^{\top}.
\end{align*}
Then, the difference between the two solutions $\tilde{\Sigma}_t:=\Sigma^2_t-\Sigma^1_t$, for all $t$, satisfies 
\begin{align}
	\tilde{\Sigma}_{t+1} = \overline{A}_t\tilde{\Sigma}_t\overline{A}^{\top}_t - \overline{A}_t\tilde{\Sigma}_tC^{\top}(\tilde{V}_t + C\tilde{\Sigma}_tC^{\top})^{-1}C\tilde{\Sigma}_t\overline{A}^{\top}_t, 
\end{align}
where $\tilde{V}_t = V+C\Sigma^1_tC^{\top}$ and $\overline{A}_t:=A-A\Sigma^1_tC^{\top}\tilde{V}_t^{-1}C$.
\end{lemma}

Next, identify $\Sigma^1_t$ with $\Sigma^*$ and $\Sigma^2_t$ with $\Sigma^*_t$ in Lemma \ref{lemma:comparison}. Then, $\tilde{V}_t = V+C\Sigma^*C^{\top} =:\tilde{V}$ and $\overline{A}_t = A-A\Sigma^*C^{\top}\tilde{V}^{-1}C =: \overline{A}$ for all $t$. Invoking Lemma \ref{lemma:comparison} leads to 
\begin{align}
	\tilde{\Sigma}_{t+1} &= \overline{A}\tilde{\Sigma}_t\overline{A}^{\top} - \overline{A}\tilde{\Sigma}_tC^{\top}(\tilde{V} + C\tilde{\Sigma}_tC^{\top})^{-1}C\tilde{\Sigma}_t\overline{A}^{\top} \label{eqn:RDE_conv_step1}\\
	&= \overline{A}\tilde{\Sigma}_t^{1/2}\big[\bI + \tilde{\Sigma}^{1/2}_tC^{\top}\tilde{V}^{-1}C\tilde{\Sigma}_t^{1/2}\big]^{-1}\tilde{\Sigma}_t^{1/2}\overline{A}^{\top} \nonumber \\
	&\hspace{-2.3em} \leq \hspace{-0.1em}\big[1\hspace{-0.15em}+\hspace{-0.15em}\lambda_{\min}(\tilde{\Sigma}^{1/2}_tC^{\top}\tilde{V}^{-1}C\tilde{\Sigma}_t^{1/2})\big]^{-1}\hspace{-0.1em}\overline{A}\tilde{\Sigma}_t\overline{A}^{\top} \hspace{-0.3em}=:\hspace{-0.1em} \mu_t \overline{A}\tilde{\Sigma}_t\overline{A}^{\top}\hspace{-0.5em},\label{eqn:RDE_conv_step2}
\end{align}
where $\tilde{\Sigma}_t^{1/2}$ denotes the unique psd square root of the psd matrix $\tilde{\Sigma}_t$, $0 < \mu_t \leq 1$ for all $t$, and $\overline{A}$ is the closed-loop matrix of the optimal infinite-horizon KF that has all its eigenvalue inside the unit circle (i.e., $\rho(\overline{A}) < 1$). Next, we use $\|\cdot\|_*$ to represent the $\Sigma^*$-induced matrix norm defined as
\begin{align*}
	\|X\|^2_* := \max_{z \neq 0} \frac{z^{\top}X^{\top}\Sigma^*Xz}{z^{\top}\Sigma^*z}.
\end{align*} 
Then, we invoke Theorem 14.4.1 of \cite{hassibi1999indefinite}, where our $\tilde{\Sigma}_t$, $\overline{A}$ and $\Sigma^*$ correspond to $P_i - P^*$, $F_p$ and $W$ in \cite{hassibi1999indefinite}, respectively. By Theorem 14.4.1 of \cite{hassibi1999indefinite} and \eqref{eqn:RDE_conv_step2}, we obtain $\|\overline{A}\|_* < 1$ and given that $\mu_t \leq 1$, 
\begin{align*}
	\|\tilde{\Sigma}_{t+1}\|_* \leq \|\overline{A}\|^2_* \cdot \|\tilde{\Sigma}_t\|_*.
\end{align*}
Therefore, the convergence rate is exponential in the sense that $\|\tilde{\Sigma}_t\|_* \leq \|\overline{A}\|_*^{2t}\cdot \|\tilde{\Sigma}_{0}\|_*$. Next, recall the condition number of a matrix $X$ is defined as $\kappa_X := \sigma_{\max}(X)/\sigma_{\min}(X)$. Then, the convergence of $\tilde{\Sigma}_t$ to the zero matrix in spectral norm can be characterized as
\begin{align*}
	\|\tilde{\Sigma}_t\| \leq \kappa_{\Sigma^*}\cdot \|\tilde{\Sigma}_t\|_* \leq \kappa_{\Sigma^*}\cdot\|\overline{A}\|_*^{2t}\cdot \|\tilde{\Sigma}_{0}\|_*.
\end{align*}
In other words, to ensure $\|\tilde{\Sigma}_N\| \leq \epsilon$, it suffices to require
\begin{align}\label{eqn:required_time}
	N \geq \frac{1}{2}\cdot \frac{\log\big(\frac{\|\tilde{\Sigma}_0\|_*\cdot \kappa_{\Sigma^*}}{\epsilon}\big)}{\log\big(\frac{1}{\|\overline{A}\|_*}\big)}.
\end{align}
Furthermore, since $(A, W^{1/2})$ is controllable, $(C, A)$ is observable, and $X_0 > \Sigma^*$, the closed-loop system at any time $t \geq 0$ is exponentially asymptotically stable such that the time-invariant (frozen) filter satisfies $\rho(A-A\Sigma^*_tC^{\top}(V+C\Sigma^*_tC^{\top})^{-1}C) < 1$ \cite{bitmead1985monotonicity, de1989monotonicity}.  Lastly, we show that the (monotonic) convergence of the filter gain to the time-invariant Kalman gain follows from the convergence of $\Sigma^*_t$ to $\Sigma^*$, which can be verified through:
\begin{align}
	L^*_t - L^* &= A\Sigma^*_tC^{\top}(V+C\Sigma^*_tC^{\top})^{-1} - A\Sigma^*C^{\top}(V+C\Sigma^*C^{\top})^{-1} \nonumber\\
 	&=A\Sigma^* C^{\top}\big[ (V + C\Sigma^*_t C^{\top})^{-1} - (V + C\Sigma^* C^{\top})^{-1}\big] \nonumber\\
 	&\hspace{1em} + A(\Sigma^*_t - \Sigma^*) C^{\top}(V + C\Sigma^*_t C^{\top})^{-1}\nonumber\\
 	&\hspace{-3em}= A\Sigma^* C^{\top}(V + C\Sigma^* C^{\top})^{-1}C(\Sigma^*-\Sigma^*_t)C^{\top}(V + C\Sigma^*_t C^{\top})^{-1} \nonumber\\
 	&\hspace{-2em}-A(\Sigma^*- \Sigma^*_t) C^{\top}(V + C\Sigma^*_t C^{\top})^{-1} \nonumber\\
 	&\hspace{-3em}=(L^*C - A)(\Sigma^* - \Sigma^*_t) C^{\top}(V + C\Sigma^*_tC^{\top})^{-1}. \label{eqn:ldiff}
\end{align}
Hence, we have $\|L^*_t - L^*\| \leq \frac{\|\overline{A}\|\cdot \|C\|}{\lambda_{\min}(V)}\cdot \|\Sigma^*_t - \Sigma^*\|$. Substituting $\epsilon$ in \eqref{eqn:required_time} with $\frac{\epsilon\cdot\lambda_{\min}(V)}{\|\overline{A}\|\cdot\|C\|}$ and identify that $\overline{A}$ is exactly $A_L^*$ completes the proof.

\subsection{Proof of Theorem \ref{theorem:KF_DP}}\label{proof:KF_DP}
To prove $\big\|[\tilde{A}_{L_{N-1}} \ \tilde{B}_{L_{N-1}}] - [A^*_L \ B^*_L]\big\| \leq \epsilon$, it suffices to bound the error between the approximated filter gain $\tilde{L}_{N-1}$ and the exact Kalman gain $L^*$ as in \eqref{eqn:kalman_gain}. First, according to Theorem \ref{lemma:finite_approximation}, we select
\begin{align}\label{eqn:N_choice}
	N = \frac{1}{2}\cdot \frac{\log\big(\frac{2\|X_0-\Sigma^*\|_*\cdot\kappa_{\Sigma^*}\cdot \|A_L^*\|\cdot\|C\|} {\epsilon\cdot\lambda_{\min}(V)}\big)}{\log\big(\frac{1}{\|A_L^*\|_*}\big)} + 1.
\end{align} 
which ensures that $L^*_{N-1}$ is stabilizing and $\|L^*_{N-1} - L^*\| \leq \epsilon/2$. Then, it remains to show that Algorithm \ref{alg:DPfilter} returns a filter $\tilde{L}_{N-1}$ such that $\|\tilde{L}_{N-1} - L^*_{N-1}\| \leq \epsilon/2$. 

Recall that the FRDE is the following forward iteration starting with $\Sigma^*_0 = X_0 > 0$: 
\begin{align}
	&\hspace{-0.4em}\Sigma^*_{t+1} = A\Sigma^*_{t}A^{\top} - A\Sigma^*_t C^{\top}(V + C\Sigma^*_t C^{\top})^{-1}C\Sigma^*_t A^{\top} + W \label{eqn:standard_RDE_appen}\\
	&\hspace{1.5em}= (A-L^*_tC)\Sigma^*_tA^{\top} + W \label{eqn:filter_riccati_finite_appen}\\
	&\hspace{1.5em}=(A-L^*_tC)\Sigma^*_t(A-L^*_tC)^{\top} + L^*_tV(L^*_t)^{\top}+W \label{eqn:filter_RDE_Lya}. 
\end{align}
 Moreover, for an arbitrary $L_t$, it holds that:
 \begin{align}
 	\Sigma_{t+1} = (A-L_tC)\Sigma_t(A-L_tC)^{\top} + L_tVL_t^{\top}+W. \label{eqn:filter_lyapunov}
 \end{align}
Furthermore, for clarity of the proof, we define/recall:
\small
\begin{align*}
 	&L^*_t\text{: Exact Kalman gain at time $t$ defined in \eqref{eqn:kalman_gain_finite}}\\
 	&\tilde{L}_t^*\text{: Optimal gain of the current cost-to-come function,}\\
 	&\hspace{1.6em} \text{absorbing errors in prior steps} \\
 	&\tilde{L}_t\text{: An approximation of $\tilde{L}_t^*$ obtained by applying \eqref{eqn:PG}}\\
 	&\delta_t:=\tilde{L}_t-\tilde{L}_t^* \text{: Policy optimization error at time $t$} \\
 	&\tilde{\Sigma}^*_{t+1} \text{: Solution generated by \eqref{eqn:filter_RDE_Lya} with $L^*_t = \tilde{L}^*_t$ and $\Sigma^*_t = \tilde{\Sigma}_t$.}
 \end{align*}
 \normalsize

We argue that $\|\tilde{L}_{N-1} - L^*_{N-1}\| \leq \epsilon/2$ can be achieved by carefully controlling $\delta_t$ for all $t$. At $t=N-1$, it holds that
\small
\begin{align*}
	\|\tilde{L}_{N-1} - L^*_{N-1}\| &\leq \|\tilde{L}^*_{N-1} - L^*_{N-1}\| + \|\delta_{N-1}\|,
\end{align*}
\normalsize
where from \eqref{eqn:ldiff} we have
\small
\begin{align*}
	\tilde{L}^*_{N-1} \hspace{-0.15em}-\hspace{-0.15em} L^*_{N-1} \hspace{-0.15em}=\hspace{-0.15em} (L^*_{N-1}\hspace{-0.1em}C \hspace{-0.1em}-\hspace{-0.1em} A)(\Sigma^*_{N-1} \hspace{-0.1em}-\hspace{-0.1em} \tilde{\Sigma}_{N-1})C^{\hspace{-0.1em}\top}\hspace{-0.1em}(V \hspace{-0.1em}+\hspace{-0.1em} C\tilde{\Sigma}_{N-1}C^{\top})^{-1}\hspace{-0.2em}.
\end{align*}
\normalsize
By Theorem \ref{lemma:finite_approximation}, $L^*_{N-1}$ is stabilizing and $\Sigma^*_{N-1} \geq \Sigma^*$ holds. Then, we require $\|\Sigma^*_{N-1} - \tilde{\Sigma}_{N-1}\| \leq \|\Sigma^*\|$ to ensure the positive definiteness of $\tilde{\Sigma}_{N-1}$. Furthermore, we can derive
\begin{align}
	\|\tilde{L}^*_{N-1} - L^*_{N-1}\| \leq \frac{\|A_{L_{N-1}}^*\|\cdot\|C\|}{\lambda_{\min}(V)}\cdot\|\Sigma^*_{N-1} - \tilde{\Sigma}_{N-1}\| \label{eqn:laststep_req}.
\end{align}
Define the helper constant
\small
\begin{align*}
	C_1:= \frac{\varphi\|C\|}{\lambda_{\min}(V)} > 0, \quad \varphi := \max_{t\in\{0, \cdots, N-1\}}\|A^*_{L_t}\|.
\end{align*}
\normalsize
Next, require $\|\delta_{N-1}\| \leq \frac{\epsilon}{4}$ and $\|\tilde{L}^*_{N-1} - L^*_{N-1}\| \leq \frac{\epsilon}{4}$ to fulfill $\|\tilde{L}_{N-1} - L^*_{N-1}\| \leq \frac{\epsilon}{2}$. By \eqref{eqn:laststep_req}, this is equivalent to requiring
\small
\begin{align}\label{eqn:sigmareq1}
	\|\Sigma^*_{N-1} - \tilde{\Sigma}_{N-1}\| \leq \min\bigg\{\|\Sigma^*\|, \frac{\epsilon}{4 C_1}\bigg\}.
\end{align}
\normalsize
Subsequently, by \eqref{eqn:filter_lyapunov}, we have
\small
\begin{align}\label{eqn:sigmalast}
	\Sigma^*_{N-1} - \tilde{\Sigma}_{N-1} = (\Sigma^*_{N-1} - \tilde{\Sigma}^*_{N-1}) + (\tilde{\Sigma}^*_{N-1} - \tilde{\Sigma}_{N-1}),
\end{align}
\normalsize
where the first difference term on the RHS of \eqref{eqn:sigmalast} is
\small
\begin{align}
\Sigma^*_{N-1} - \tilde{\Sigma}^*_{N-1} &= (A-L^*_{N-2}C)\Sigma^*_{N-2}A^{\top} - (A-\tilde{L}^*_{N-2}C)\tilde{\Sigma}_{N-2}A^{\top} \nonumber\\
&\hspace{-6.5em}= (A\hspace{-0.15em}-\hspace{-0.15em}L^*_{N\hspace{-0.1em}-\hspace{-0.1em}2}C)(\Sigma^*_{N\hspace{-0.1em}-\hspace{-0.1em}2} \hspace{-0.15em}-\hspace{-0.15em}\tilde{\Sigma}_{N-2})A^{\hspace{-0.15em}\top} \hspace{-0.3em}+ (\tilde{L}^*_{N\hspace{-0.1em}-\hspace{-0.1em}2} \hspace{-0.15em}-\hspace{-0.15em} L^*_{N\hspace{-0.1em}-\hspace{-0.1em}2})C\tilde{\Sigma}_{N-2}A^{\hspace{-0.15em}\top}\hspace{-0.3em}.\label{eqn:sigmalast2}
\end{align}
\normalsize
Moreover, the second term on the RHS of \eqref{eqn:sigmalast} is 
\small
\begin{align}
	&\tilde{\Sigma}^*_{N-1} \hspace{-0.15em}-\hspace{-0.05em} \tilde{\Sigma}_{N-1} \hspace{-0.1em}=\hspace{-0.1em} (A\hspace{-0.1em}-\hspace{-0.1em}\tilde{L}^*_{N-2}C)\tilde{\Sigma}_{N-2}(A\hspace{-0.1em}-\hspace{-0.1em}\tilde{L}^*_{N-2}C)^{\hspace{-0.1em}\top} \hspace{-0.15em}+\hspace{-0.05em} \tilde{L}^*_{N-2}V(\tilde{L}^*_{N-2})^{\hspace{-0.1em}\top} \nonumber\\
	&\hspace{1em} - (A-\tilde{L}_{N-2}C)\tilde{\Sigma}_{N-2}(A-\tilde{L}_{N-2}C)^{\top} - \tilde{L}_{N-2}V(\tilde{L}_{N-2})^{\top} \nonumber\\
	&= A\tilde{\Sigma}_{N-2}C^{\top}\tilde{L}_{N-2}^{\top} - A\tilde{\Sigma}_{N-2}C^{\top}(\tilde{L}^*_{N-2})^{\top} + \tilde{L}_{N-2}C\tilde{\Sigma}_{N-2}A^{\top} \nonumber\\
	&\hspace{1em}+ \tilde{L}^*_{N-2}(C\tilde{\Sigma}_{N-2}C^{\top} + V)(\tilde{L}^*_{N-2})^{\top} - \tilde{L}_{N-2}(C\tilde{\Sigma}_{N-2}C^{\top} + V)\tilde{L}_{N-2}^{\top} \nonumber\\
	&\hspace{1em}- \tilde{L}^*_{N-2}C\tilde{\Sigma}_{N-2}A^{\top} \nonumber\\
	&=\tilde{L}^*_{N-2}(C\tilde{\Sigma}_{N-2}C^{\top} + V)(\tilde{L}^*_{N-2})^{\top}- A\tilde{\Sigma}_{N-2}C^{\top}(\tilde{L}^*_{N-2})^{\top} \nonumber\\
	&\hspace{1em} - \tilde{L}^*_{N-2}C\tilde{\Sigma}_{N-2}A^{\top}  + A\tilde{\Sigma}_{N-2}C^{\top}(C\tilde{\Sigma}_{N-2}C^{\top}+V)^{-1}C\tilde{\Sigma}_{N-2}A^{\top}\nonumber\\
	&\hspace{1em} -(\tilde{L}_{N-2} - A\tilde{\Sigma}_{N-2}C^{\top}(C\tilde{\Sigma}_{N-2}C^{\top}+V)^{-1})(C\tilde{\Sigma}_{N-2}C^{\top}+V)\nonumber\\
	&\hspace{2em}\cdot(\tilde{L}_{N-2} - A\tilde{\Sigma}_{N-2}C^{\top}(C\tilde{\Sigma}_{N-2}C^{\top}+V)^{-1})^{\top} \label{eqn:sigmadiff1}\\
	&= (A\tilde{\Sigma}_{N-2}C^{\top}(C\tilde{\Sigma}_{N-2}C^{\top}+V)^{-1} - \tilde{L}^*_{N-2})(C\tilde{\Sigma}_{N-2}C^{\top}+V)\nonumber\\
	&\hspace{2em}\cdot(A\tilde{\Sigma}_{N-2}C^{\top}(C\tilde{\Sigma}_{N-2}C^{\top}+V)^{-1}-\tilde{L}^*_{N-2})^{\top}\nonumber\\
	&\hspace{1em} -(\tilde{L}_{N-2} - A\tilde{\Sigma}_{N-2}C^{\top}(C\tilde{\Sigma}_{N-2}C^{\top}+V)^{-1})(C\tilde{\Sigma}_{N-2}C^{\top}+V)\nonumber\\
	&\hspace{2em}\cdot(\tilde{L}_{N-2} - A\tilde{\Sigma}_{N-2}C^{\top}(C\tilde{\Sigma}_{N-2}C^{\top}+V)^{-1})^{\top}, \label{eqn:sigmadiff2}
\end{align}
\normalsize
where \eqref{eqn:sigmadiff1} is due to completion of squares. Substituting $\tilde{L}^*_{N-2} = A\tilde{\Sigma}_{N-2}C^{\top}(C\tilde{\Sigma}_{N-2}C^{\top}+V)^{-1}$ into \eqref{eqn:sigmadiff2} leads to
\begin{align}\label{eqn:sigmadiff3}
	\tilde{\Sigma}^*_{N-1} - \tilde{\Sigma}_{N-1} = -\delta_{N-2}(C\tilde{\Sigma}_{N-2}C^{\top}+V)\delta_{N-2}^{\top}.
\end{align}
Thus, combining \eqref{eqn:sigmalast}, \eqref{eqn:sigmalast2}, and \eqref{eqn:sigmadiff3} yields
\small
\begin{align}
	&\hspace{1em}\|\Sigma^*_{N-1} - \tilde{\Sigma}_{N-1}\| \nonumber\\
	&\leq \|\Sigma^*_{N-2} -\tilde{\Sigma}_{N-2}\|\cdot\varphi\|A\|  + \|\tilde{L}^*_{N-2} - L^*_{N-2}\|\cdot\|C\|\cdot\|\tilde{\Sigma}_{N-2}\|\cdot\|A\| \nonumber\\
	&\hspace{1em}+ \|\delta_{N-2}\|^2\|C\tilde{\Sigma}_{N-2}C^{\top}+V\| \nonumber \\
	&\leq\|A\|\cdot[\varphi+C_1\cdot\|C\|\cdot\|\tilde{\Sigma}_{N-2}\|]\cdot\|\Sigma^*_{N-2} - \tilde{\Sigma}_{N-2}\| \nonumber\\
	&\hspace{1em}+ \|\delta_{N-2}\|^2\|C\tilde{\Sigma}_{N-2}C^{\top}+V\| \label{eqn:sigma_contraction},
\end{align}
\normalsize
where the last inequality follows from \eqref{eqn:laststep_req}. Now, require
\small
\begin{align}
	\|\Sigma^*_{N-2} - \tilde{\Sigma}_{N-2}\| &\leq \min\bigg\{\|\Sigma^*\|, \frac{\|\Sigma^*\|}{C_2}, \frac{\epsilon}{4 C_1C_2}\cdot\bigg\} \label{eqn:sigmareq2}\\
	\|\delta_{N-2}\| &\leq  \min\bigg\{\sqrt{\frac{\|\Sigma^*\|}{C_3}}, \frac{1}{2}\sqrt{\frac{\epsilon}{C_1C_3}}\bigg\}\label{eqn:L_req1},
\end{align}
\normalsize
where $C_2$ and $C_3$ are positive constants defined as
\begin{align*}
	C_2 &:= 2\|A\|\cdot\big[\varphi+C_1\cdot\|C\|\cdot\big(\|X_0\| + \|\Sigma^*\|\big)\big] >0 \\
	C_3 &:= 2\big[\|V\| + \|C\|^2\big(\|X_0\|+ \|\Sigma^*\|\big)\big] > 0.
\end{align*}
Then, conditions \eqref{eqn:sigmareq2} and \eqref{eqn:L_req1} are sufficient for \eqref{eqn:sigmareq1} (and thus for $\|\tilde{L}_{N-1} - L^*_{N-1}\| \leq \epsilon/2$) to hold. Subsequently, we can propagate the required accuracies in \eqref{eqn:sigmareq2} and \eqref{eqn:L_req1} backward in time. Specifically,   we iteratively apply the arguments in \eqref{eqn:sigma_contraction} (i.e., by plugging quantities with subscript $t$ into the LHS of \eqref{eqn:sigma_contraction} and plugging quantities with subscript $t-1$ into the RHS of \eqref{eqn:sigma_contraction}) to obtain  the result that if at all $t \in \{1, \cdots, N-2\}$, we require
\small
\begin{align}
	&\|\Sigma^*_{t} - \tilde{\Sigma}_{t}\| \leq \min\bigg\{\|\Sigma^*\|, \frac{\|\Sigma^*\|}{C_2^{N-t-1}}, \frac{\epsilon}{4 C_1C_2^{N-t-1}}\bigg\} \label{eqn:sigmareq_allt}\\
	&\hspace{-0.5em}\|\delta_{t}\| \leq  \min\bigg\{\sqrt{\frac{\|\Sigma^*\|}{C_3}}, \sqrt{\frac{\|\Sigma^*\|}{C_2^{N-t-2}C_3}}, \frac{1}{2}\sqrt{\frac{\epsilon}{C_1C_2^{N-t-2}C_3}}\bigg\} \nonumber,
\end{align}
\normalsize
then \eqref{eqn:sigmareq2} holds true and therefore \eqref{eqn:sigmareq1} is satisfied. 

We now compute the required accuracy for $\delta_0$. As illustrated in Figure \ref{fig:proof_sketch}, we have $\Sigma^*_1 = \tilde{\Sigma}^*_1$ because 
\small
\begin{align*}
	\tilde{\Sigma}^*_1 &= (A-\tilde{L}^*_0C)\Sigma^*_0(A-\tilde{L}^*_0C)^{\top} + \tilde{L}^*_0V(\tilde{L}^*_0)^{\top}+W \\
	&= (A-L^*_0C)\Sigma^*_0(A-L^*_0C)^{\top} + L^*_0V(L^*_0)^{\top}+W = \Sigma^*_1,
\end{align*}
\normalsize
where the second equality is due to $\tilde{L}^*_0 = L^*_0$ since there are no prior computational errors yet at $t=0$. By \eqref{eqn:sigma_contraction}, the distance between $\Sigma^*_1$ and $\tilde{\Sigma}_1$ can be bounded as 
\small
\begin{align*}
	\|\Sigma^*_1 - \tilde{\Sigma}_1\| = \|\tilde{\Sigma}^*_1 - \tilde{\Sigma}_1\| \leq \|\delta_{0}\|^2\cdot C_3.
\end{align*}
\normalsize
To fulfill the requirement \eqref{eqn:sigmareq_allt} for $t=1$, which is
\small
\begin{align*}
	\|\Sigma^*_{1} - \tilde{\Sigma}_{1}\| \leq \min\bigg\{\|\Sigma^*\|, \frac{\|\Sigma^*\|}{C_2^{N-2}}, \frac{\epsilon}{4 C_1C_2^{N-2}}\bigg\},
\end{align*}
\normalsize
it suffices to let
\small
\begin{align}
	\|\delta_0\| \leq \min\bigg\{\sqrt{\frac{\|\Sigma^*\|}{C_3}}, \sqrt{\frac{\|\Sigma^*\|}{C_2^{N-2}C_3}}, \frac{1}{2}\sqrt{\frac{\epsilon}{C_1C_2^{N-2}C_3}}\bigg\}. \label{eqn:delta0_req}
\end{align}
\normalsize

Finally, we analyze the worst-case complexity of the proposed algorithm by computing, at the most stringent case, the required size of $\|\delta_t\|$. When $C_2 \leq 1$, the most stringent dependence of $\|\delta_t\|$ on $\epsilon$ happens at $t=N-1$, which is of the order $\cO(\epsilon)$, and the dependences on system parameters (through the dependence on constants $C_1, C_2$ and $C_3$) are polynomial. We argue that if $C_2 > 1$, then the requirement on $\|\delta_{N-1}\|$ is still the most stringent one. This is because $\|\delta_0\| \leq \|\delta_t\|$ for all $t \in \{1, \cdots, N-2\}$ and by \eqref{eqn:delta0_req}, we have
\small
\begin{align}\label{eqn:delta0}
	\|\delta_0\| \sim \cO\Big(\sqrt{\frac{\epsilon}{C_2^{N-2}}}\Big).
\end{align}
\normalsize
Since we require $N$ to satisfy \eqref{eqn:N_choice}, the dependence of $\|\delta_0\|$ on $\epsilon$ in \eqref{eqn:delta0} becomes $\|\delta_0\| \sim \cO(\epsilon^{\frac{3}{4}})$, which is milder than that of $\|\delta_{N-1}\|$. Therefore, it suffices to require the most stringent error bound for all $t$, which is $\|\delta_t\| \sim \cO(\epsilon)$, to reach the $\epsilon$-neighborhood of the infinite-horizon KF gain. Lastly, for $\tilde{A}_{L_{N-1}}$ to be stabilizing, it suffices to select a sufficiently small $\epsilon$ such that
\small
\begin{align*}
	\epsilon < 1 - \|A_L^*\|_* \Longrightarrow \|\tilde{A}_L\|_* < 1.
\end{align*}
\normalsize
This completes the proof.

\end{document}